\documentstyle[12pt]{article}
\title{\LARGE\bf On the weak Freese-Nation property of complete Boolean
algebras}
\author{Saka\'e Fuchino, Stefan Geschke, Saharon Shelah, Lajos Soukup%
\ifcommented
\bigskip\bigskip\\
\fbox{\,\parbox{7.6cm}{\footnotesize\tt\mbox{}\hfill
This is still a preliminary version.\hfill\mbox{}\\
\mbox{}\hfill Any comments are appreciated.\hfill\mbox{}
}\,}
\fi
}
\newif\iftesting
\newif\ifcommented
\commentedtrue
\font\teneufb=eufb10
\font\seveneufb=eufb7
\font\fiveeufb=eufb5
\newfam\eufbfam
\textfont\eufbfam=\teneufb
\scriptfont\eufbfam=\seveneufb
\scriptscriptfont\eufbfam=\fiveeufb

\font\tenmsbm=msbm10
\font\sevenmsbm=msbm7
\font\fivemsbm=msbm5
\newfam\msbmfam
\textfont\msbmfam=\tenmsbm
\scriptfont\msbmfam=\sevenmsbm
\scriptscriptfont\msbmfam=\fivemsbm
\def\bbd#1{{\fam\msbmfam\relax#1}}
\iftesting
\let\Label\label%
\def\label#1{\mbox{}\marginpar{{\tiny #1}}\Label{#1}\ignorespaces}%
\fi

\ifcommented
\fi
{\ifcommented\end{footnotesize}\medskip\\\fi}

\newtheorem{Thm}{{\bf Theorem}}[section]
\newtheorem{Cor}[Thm]{{\bf Corollary}}
\newtheorem{Prop}[Thm]{{\bf Proposition}}
\newtheorem{Lemma}[Thm]{{\bf Lemma}}

\newtheorem{Claim}{{\bf Claim}}[Thm]
\newtheorem{Problem}{{\bf Problem}}

\newtheorem{qs}{Question}

\newcommand{\Thmof}[1]{{Theorem \ref{#1}}}
\newcommand{\Corof}[1]{{Corollary \ref{#1}}}
\newcommand{\Propof}[1]{{Proposition \ref{#1}}}
\newcommand{\Lemmaof}[1]{{Lemma \ref{#1}}}

\newcommand{\Claimof}[1]{{Claim \ref{#1}}}

\newcommand{\prf}{{\bf Proof\ \ }\ignorespaces}

\newcommand{\prfofClaim}{\raisebox{-.4ex}{\Large $\vdash$\ \ }}
\newsavebox{\qedbox}\sbox{\qedbox}{
S.F.)
{\unitlength=0.07mm \begin{picture}(40,60)
\put(0,0){\framebox(30,44)[cc]{}}
\put(30,-7){\rule{7\unitlength}{44\unitlength}}
\put(10,-7){\rule{27\unitlength}{7\unitlength}}
\end{picture}}}
\newcommand{\qed}{\mbox{}\hfill\usebox{\qedbox}}
\newcommand{\smallqed}%
{\mbox{}\smallskip\hfill\raisebox{-.4ex}{\Large $\dashv$}\\}
\newcommand{\qedof}[1]%
{\mbox{} \hspace*{\fill}{\usebox{\qedbox}{~(#1)}}%
\mbox{}}
\newcommand{\Qedof}[1]%
{\mbox{} \hspace*{\fill}{\usebox{\qedbox}%
{~(#1~\number\theThm)}}}
\newcommand{\qedofThm}{\Qedof{Theorem}}
\newcommand{\qedofCor}{\Qedof{Corollary}}
\newcommand{\qedofProp}{\Qedof{Proposition}}

\newcommand{\qedskip}{\medskip}
\newcommand{\qedofClaim}%
{\mbox{}\hfill\raisebox{-.4ex}{\Large $\dashv$ }\nolinebreak%
\mbox{~(Claim~\number\theClaim)}}

\newcommand{\st}{such that}
\newcommand{\wrt}{with respect to}

\newcommand{\wolog}{without loss of generality}
\newcommand{\Wolog}{Without loss of generality}
\newcommand{\Tfae}{The following are equivalent}
\newcommand{\tfae}{the following are equivalent}
\newcommand{\po}{partial ordering}
\newcommand{\pos}{partial orderings}
\newcommand{\tenten}{\ifmmode,\ldots,\else\mbox{\rm,\ldots,} \fi}

\newcommand{\Ba}{Boolean algebra}
\newcommand{\Bas}{Boolean algebras}
\newcommand{\cBa}{complete Boolean algebra}
\newcommand{\cBas}{complete Boolean algebras}
\newcommand{\wfn}{weak Freese-Nation property}
\newcommand{\wfnm}{weak Freese-Nation mapping}
\newcommand{\mild}{neat}
\newcommand{\Implies}{\,$\Rightarrow$\,}

\newenvironment{assertion}[1]{\begin{trivlist}
\newbox\assertbox
\dimen255=\textwidth
\setbox\assertbox=\hbox{\hspace*{\parindent}{#1}\hspace{\labelsep}}
\advance\dimen255 by -1\wd\assertbox
\advance\dimen255 by -1ex
\item[]\unhbox\assertbox\hfill
\begin{minipage}[t]{\dimen255}}%
{\end{minipage}\end{trivlist}}
{\hfill\mbox{}\end{trivlist}}
\def\assert#1{\makebox[5ex][l]{\rm (#1)}\ignorespaces}
\def\lassert#1{\llap{\makebox[5ex][l]{\rm (#1)}}\ignorespaces}
\def\assertof#1{{\rm (#1)}}

\newenvironment{markedformula}[1]%
        {\begin{trivlist}\item[\,\,#1]\mbox{}\hfill$\displaystyle}%
        {$\hfill\mbox{}\end{trivlist}}
\newcommand{\setof}[2]{\{#1\,:\,#2\}}
\newcommand{\ssetof}[1]{\{#1\}}
\newcommand{\seqof}[2]{\langle#1\,:\,#2\rangle}

\newcommand{\psetof}[1]{{\cal P}\/(#1)}
\newcommand{\psof}[1]{{\cal P}\/(#1)}
\newcommand{\cardof}[1]{\mathopen{|\,}#1\mathclose{\,|}}
\newcommand{\mapping}[3]{#1:#2\rightarrow #3}
\newcommand{\pfeil}{\mathrel{{\rightarrow\,}\llap{$\rightarrow$}}}
\newcommand{\fnsp}[2]{\mbox{}^{#1^{\mbox{}\!}}#2}
\newcommand{\bairesp}{\fnsp{\omega}{\omega}}

\newcommand{\restr}%
{{\hspace{0.1ex}|\hspace{-0.02ex}{\grave{}}\hspace{0.8ex}}}

\renewcommand{\Re}{\bbd{R}}
\newcommand{\rationals}{\bbd{Q}}

\newcommand{\regembed}%
        {\mathrel{{<}\llap{\raisebox{0.2ex}{$\scriptstyle\,\circ$}}}}
\newcommand{\forces}[2]{\,\|\hspace{-.35ex}\mbox{\sf--}_{\,#1\,}%
\mbox{\rm``}\,#2\,\mbox{\rm''}}

\newcommand{\dotleq}{\mathrel{\dot{\leq}}}
\newcommand{\lessnoneq}%
{\mathrel{\raisebox{-0.8ex}{$\stackrel{<}{\scriptstyle\,\not=\,}$}}}

\newcommand{\xmbox}[1]{ $\relax{\rm #1}\relax$ }
\newcommand{\dom}{\mathop{\rm dom}}

\newcommand{\cof}{\mathop{\rm cof}}
\newcommand{\cf}{\mathop{\rm cf}}
\newcommand{\Fn}{\mathop{\rm Fn}}

\newcommand{\ZFC}{{\rm ZFC}}
\newcommand{\CH}{{\rm CH}}
\newcommand{\GCH}{{\rm GCH}}

\newcommand{\calH}{{\cal H}}

\newcommand{\calS}{{\cal S}}

\newcommand{\dotd}{{\dot d}}

\newcommand{\dotg}{{\dot g}}
\newcommand{\doth}{{\dot h}}
\newcommand{\dotx}{{\dot x}}

\newcommand{\dotF}{{\dot F}}

\newcommand{\dotN}{{\dot N}}
\newcommand{\dotQ}{{\dot Q}}
\newcommand{\dotS}{{\dot S}}

\newcommand{\dotvarphi}{{\dot\varphi}}
\newcommand{\doteta}{{\dot\eta}}
\begin{document}
\maketitle
\begin{abstract}
The following results are proved:

(a)
{\it In a model obtained by adding $\aleph_2$ Cohen reals, there is
always a c.c.c.\  \cBa\  without the \wfn. }

(b) {\it Modulo the consistency strength of a supercompact cardinal, the
existence of a c.c.c.\  \cBa\  without the \wfn\  is consistent with \GCH. 
}

(c) {\it If a weak form of $\Box_\mu$ and
$\cof([\mu]^{\aleph_0},{\subseteq})=\mu^+$ hold for eavch
$\mu>\cf(\mu)=\omega$, then
the \wfn\  of
$(\psof{\omega},{\subseteq})$ is equivalent to the \wfn\  of any of
$\bbd{C}(\kappa)$ or $\bbd{R}(\kappa)$ for uncountable $\kappa$. }

(d) {\it Modulo consistency of
$(\aleph_{\omega+1},\aleph_\omega)\pfeil(\aleph_1,\aleph_0)$, it is
consistent with \GCH\  that $\bbd{C}(\aleph_\omega)$ does not have the
\wfn\  and hence
the assertion in
(c) does not hold, and also that adding $\aleph_\omega$ Cohen reals 
destroys
the \wfn\  of $(\psof{\omega},{\subseteq})$. }\medskip

These results solve all of the problems listed in Fuchino-Soukup
\cite{fuchino-soukup} and some other problems posed by S.\  Geschke.
\def\thefootnote{\fnsymbol{footnote}}
\setcounter{footnote}{1}
\footnotetext{The first author was
partially supported by Grant-in-Aid for
Scientific Research (C) No.\  10640099 of the Ministry of Education,
Science, Sports and Culture, Japan.

Some of the results here, in
particular earlier versions of the results in Section \ref{square},
were included in the second
author's Ph.D.\  thesis \cite{geschke}.

This paper is [FGShS:712] of the third author's publications list.
His research is supported by ``The Israel Science Foundation''.

The fourth author was partially
supported by Grant-in-Aid for JSPS Fellows No.\  98259 of the Ministry of
Education, Science, Sports and Culture, Japan, and by Hungarian National
Foundation for Scientific Research grant no.\  25745.

Section \ref{chang} of the present paper was worked out during the
conference "Algebra and Discrete Mathematics" in
Hattingen, Germany 1999 which was organized
by the third author and attended by all of the other authors.

The final version of the paper was then prepared during the
Workshop on Set-Theoretical Topology 1999 at Erd\H os Center, Budapest
Hungary 1999.
}
\end{abstract}
\section{Introduction}
A quasi-ordering $(P,\leq)$ is said to have {\it the \wfn\/} if there
is a mapping $\mapping{f}{P}{[P]^{\leq\aleph_0}}$ \st:
\begin{assertion}{} For any $p$, $q\in P$ with $p\leq q$ there is
$r\in f(p)\cap f(q)$ \st\    $p\leq r\leq q$.
\end{assertion}
A mapping $f$ as above is called a {\it\wfnm} on $P$.

The \wfn\    was introduced in Chapter 4 of \cite{heindorf-shapiro} as a
weakening of a
notion of almost freeness of \Ba s. The property was further studied
in \cite{FuKoSh} and \cite{fuchino-soukup}.

In \cite{FuKoSh}, it is shown that $(\psof{\omega_1},\subseteq)$ does not
have the \wfn. If a \cBa\    $B$ does
not have the c.c.c., then $(\psof{\omega_1},\subseteq)$ can be completely
embedded into $B$. Hence, in this case, $B$ can not have the \wfn.

It is easily seen that every quasi-ordering of cardinality
$\leq\aleph_1$ has the \wfn\,(see e.g.\    \cite{FuKoSh}). It follows that,
under \CH,
$(\psof{\omega},{\subseteq})$ has the \wfn.

To simplify the formulation of some of the results below,
let us say that  a model of set-theory is {\em \mild} if
$\Box_{{\mu}}$ holds --- what is actually needed in the following is merely 
a
very weak variant of $\Box_\mu$ introduced in \cite{fuchino-soukup} (see
before \Propof{astastast}) --- and 
$\cof([{\mu}]^{\omega},\subseteq)={\mu}^+$
for each ${\mu}>\cf({\mu})={\omega}$.

In \cite{FuKoSh} and on \cite{fuchino-soukup}, it was shown that  if CH 
holds,
then every c.c.c.
\cBa\    of size $<\aleph_{\omega}$
has the \wfn;  and in a \mild\  model,
CH implies that every c.c.c.\   \cBa\   has the \wfn.
However, the following questions remained unanswered in 
\cite{fuchino-soukup}:

\begin{qs}\label{q:eq}{\rm(\cite[Problem 5]{fuchino-soukup})}
Are the following equivalent?\\
\assert{i} $(\psof{\omega},{\subseteq})$ has the
\wfn;\\
\assert{ii} every c.c.c.\   \cBa\   has the \wfn.
\end{qs}

\begin{qs}\label{q:gch}{\rm(\cite[Problem 2]{fuchino-soukup})}
Does $\ZFC+\GCH$ imply that every  c.c.c.\    \cBa\    has  the
\wfn\,?
\end{qs}

We give negative answers here:
see \Corof{Cor2} for question \ref{q:eq} and
\Thmof{Cor3} for question \ref{q:gch}.
By the result in \cite{fuchino-soukup} already mentioned above, we need
consistency strength of some large cardinal to give a negative  answer
for Question \ref{q:gch}.
 Indeed, the ground model $V$ in the negative solution of this
problem is obtained by starting from a model of
\ZFC\    with a supercompact cardinal.

In \cite{FuKoSh} it was shown that  if CH holds, then
adding less than $\aleph_{\omega}$ many Cohen
reals preserves the \wfn\    of
$(\psof{\omega},{\subseteq})$.
By \cite{fuchino-soukup}, in the generic extension obtained by adding any
number of Cohen reals to a \mild\    model satisfying CH, 
 not only $(\psof{\omega},{\subseteq})$ but  every tame
c.c.c.\    \cBa\    has the \wfn. Here, letting $P=\Fn(\tau,2)$ ($=$ the 
standard
p.o.\    for adding $\tau$ Cohen reals), a \Ba\    in a $P$-generic
extension is said to be {\it tame}, if there is a
$P$-name $\dotleq$ of partial ordering of $B$ and a mapping
$\mapping{t}{B}{[\tau]^{\aleph_0}}$ in $V$ \st, for every $p\in P$ and $x$,
$y\in B$, if $p\forces{P}{x\dotleq y}$, then
$p\restr(t(x)\cup t(y))\forces{P}{x\dotleq y}$ (we assume here \wolog\
that $B$ is chosen so that its underlying set is a ground model set).

These results suggest the following questions posed in
\cite{fuchino-soukup}:

\begin{qs}\label{q:poc}
\label{problem3}{\rm(\cite[Problem  3]{fuchino-soukup})}
Assume that $V[G]$ is a model obtained by adding Cohen reals to a model
of $\ZFC +\CH$. Is it true that $\psof{\omega}$ has the
\wfn\  in $V[G]$\,?
\end{qs}
\begin{qs}\label{q:cot}
{\rm(\cite[Problem 4]{fuchino-soukup})}
Assume that $V[G]$ is a model obtained by adding $\aleph_2$
Cohen reals to a model
of $\ZFC +\CH$. Is it true that every
 c.c.c.\    \cBa\    (not just the tame ones) has the \wfn\    in $V[G]$?
\end{qs}

The results of the present paper answer these questions in the
 negative: see
 \Thmof{chang-pomega} for question \ref{q:poc} and
\Corof{Cor1} for question
\ref{q:cot}.

By the result in \cite{fuchino-soukup} already mentioned above, we need
consistency strength of some large cardinal to give a negative solution of
Question \ref{q:poc}.
Indeed, the ground model $V$ in the negative solution of this
problem given in \Thmof{chang-pomega}
is obtained by starting from a model of \ZFC\  with a large cardinal
slightly stronger than a huge cardinal.

After the negative solution of Problem 5, the following question still
remains:
\begin{Problem} For which \Ba\    $B$, the \wfn\    of $B$ is equivalent 
with
the \wfn\    of $(\psof{\omega},{\subseteq})$\,?
\end{Problem}

The following easy lemma is already a
result  in this direction.
\begin{Lemma}\label{equivalence}
\Tfae:\medskip\\
\assert{$a$} $(\psof{\omega},{\subseteq})$ has the \wfn;\\
\assert{$b$} $(\psof{\omega}, \subseteq^*)$ has the \wfn;\\
\assert{$c$} $(\psof{\omega}/fin, \subseteq^*)$ has the \wfn;\\
\assert{$d$} $(\bairesp,\leq)$ has the \wfn;\\
\assert{$e$} $(\bairesp,\leq^*)$ has the \wfn.\\
\assert{$f$} $({}^{\omega}\Re,\leq)$ has the \wfn.\qed
\end{Lemma}

Koppelberg \cite{koppelberg} pointed out that the \wfn\    of Cohen algebra
$\bbd{C}(\omega)$ is equivalent to
the \wfn\    of $(\psof{\omega},{\subseteq})$. In the
present paper, we show that it is
also equivalent to the \wfn\    of the measure algebra
${\bbd R}(\omega)$ (\Propof{measure-alg}) and
more over,
in a mild model
also with \wfn\
of $\bbd{C}(\kappa)$ and/or $\bbd{R}(\kappa)$ for any
$\kappa\geq\aleph_0$ (\Corof{equivalence2}). Here, we denote with
$\bbd{C}(\kappa)$ and $\bbd{R}(\kappa)$ the c.c.c.\    \cBas\
$Borel(\fnsp{\kappa}{2})/meager(\fnsp{\kappa}{2})$ and
$Borel(\fnsp{\kappa}{2})/null(\fnsp{\kappa}{2})$ respectively.
We show that
some extra set-theoretic assumption are
really necessary in \Corof{equivalence2} by constructing a model of
\GCH\    and the negation of \wfn\    for $\bbd{C}(\aleph_\omega)$ starting
from a model of \GCH\    and Chang's conjecture for $\aleph_{\omega}$.


Assume that $\langle P_{\alpha}, \dotQ_{\alpha}:{\alpha}<{\omega}_2
\rangle$
is a finite support iteration such that
forcing with $\dotQ_{\alpha}$ just adds a real to $V^{P_{\alpha}}$.
Then, as S. Geschke proved in  \cite{geschke} , if
this iteration preserves the \wfn\  of $\psof{\omega}$ then
for all but ${\omega}_1$ many ${\alpha}$ the partially ordered set 
$\dotQ_{\alpha}$
just adds one Cohen real. But
by \Corof{Cor1}, in  any model obtained by adding $\geq\aleph_2$
Cohen reals, there is a c.c.c.\   \cBa\   $B$ without the \wfn.
So there is no easy way to blow up the continuum and to preserve the \wfn\
of all c.c.c \cBas. Thus the following question
seems to be a reasonable one:

\begin{Problem}
Does \CH\    follow from the assumption that
every c.c.c.\    \cBa\    has the \wfn?
\end{Problem}

If  ${\bf b}>\aleph_1$ or if there is an $\aleph_2$-Luzin-gap, then
  $(\psof{\omega},{\subseteq})$ does not have
the \wfn\    (see \cite{FuKoSh} and \cite{fuchino-soukup}).
The following question ({\rm(\cite[Problem 1]{fuchino-soukup})})
was raised against this background: \medskip

{\em
Suppose that $\psof{\omega}$ does not have any increasing chain of
length $\geq\omega_2$ \wrt\    $\subseteq^*$ and that there is no
$\aleph_2$-Luzin gap. Does it follow that $\psof{\omega}$ has the
\wfn\,?} \medskip

This problem  can be solved negatively using results from
\cite{brendle-fuchino-soukup} and \cite{geschke}:
Let $V$ be a model of \CH\
and $V[G]$ its generic extension by adding many random reals side by side.
Using results from
\cite{brendle-fuchino-soukup} we see that in $V[G]$, there are neither
increasing $\omega_2$ chain in $\psof{\omega}$ \wrt\    $\subseteq^*$ nor
$\aleph_2$-Luzin gap. On the other hand S.\    Geschke \cite{geschke} 
showed
that in $V[G]$ $(\psof{\omega},\subseteq)$ does not have the \wfn.

Consequences of the
\wfn\    of $(\psof{\omega},\subseteq)$ were studied in \cite{koppelberg}
and \cite{fuchino-geschke-soukup}. In the latter
paper it was shown that a set-theoretic universe with the \wfn\
of $(\psof{\omega},\subseteq)$ looks quite similar to a Cohen model. In
particular, under the \wfn\    of $(\psof{\omega},\subseteq)$, all cardinal
invariants which appear in \cite{blass} take the same value as in a
Cohen model with the same size of $2^{\aleph_0}$.

\begin{Problem}
Find a combinatorial $(\mathbf \Pi^1_1)$ characterization of
\wfn\    of $\psof{\omega}$.
\end{Problem}

The \wfn\    of a quasi-ordering $(P,\leq)$ is actually a property of the
corresponding partial ordering $(\overline{P},\leq)$ obtained as the
quotient structure of $(P,\leq)$ \wrt\    the equivalence relation
``$x\leq y\land y\leq x$'': $(P,\leq)$ has the \wfn\    if and only if
$(\overline{P},\leq)$ does.

The following criteria of the \wfn\    are used in the later sections.
A \po\    $Q$ is said to be a {\it retract} of a \po\    $P$ if there are 
order
preserving mappings $\mapping{i}{Q}{P}$ and
$\mapping{j}{P}{Q}$ \st\    $j\circ i=id_Q$.
Note that if $P$ and $Q$ are \cBas\    and there is a strictly 
order-preserving
embedding $f$ of $Q$ into $P$ (i.e.\    $f$ preserves ordering and
incomparability) then we can always find order preserving
$\mapping{g}{P}{Q}$ \st\    $g\circ f=id_Q$: simply define $g$ by
$g(p)=\sum\setof{q\in Q}{f(q)\leq p}$ for $p\in P$.

$Q$ is said to be a {\it $\sigma$-subordering of\/} $P$ (notation:
$Q\leq_\sigma P$) if, for
every $p\in P$, $Q\restr p=\setof{q\in Q}{q\leq p}$ has a countable
cofinal subset and $Q\uparrow p=\setof{q\in Q}{q\geq p}$ has a countable
coinitial subset.
Note that if $C$ is a complete subalgebra of a \cBa\    $B$ (notation:
$C\leq_c B$) or a countable union of complete subalgebras of $B$, then it
follows that $C\leq_\sigma B$.

\begin{Prop}\mbox{}\assert{a} {\rm (Lemma 2.7 in \cite{FuKoSh})}
\label{criterions}
If $Q$ is a retract of $P$ and $P$ has the \wfn\    then $Q$ has the
\wfn.\smallskip\\
\assert{b} {\rm (Lemma 2.3 \assertof{a} in \cite{FuKoSh})}
If $Q\leq_\sigma P$ and $P$ has the \wfn, then $Q$ also has the
\wfn.\smallskip\\
\assert{c} {\rm (Lemma 2.6 in \cite{FuKoSh})}
If $P_\alpha$, $\alpha<\delta$ is an increasing sequence of \pos\    with 
the
\wfn\    \st\
$P_\alpha\leq_\sigma P_{\alpha+1}$ for every $\alpha<\delta$ and
$P_{\gamma}=\bigcup_{\alpha<\gamma}P_\alpha$ for all $\gamma<\delta$ with
$\cf(\gamma)>\omega$, then $P=\bigcup_{\alpha<\delta}P_\alpha$ also has
the \wfn. \qed
\end{Prop}
\section{$P_\calS$ and $B_\calS$}
In this section we introduce a construction of \pos\  $P_\calS$ and \Bas\
$B_\calS$ which will be used in Sections \ref{cohen-reals} and \ref{GCH}.
For $S\subseteq \kappa$ and an indexed family
$\calS=\seqof{S_\alpha}{\alpha\in S}$ of subsets of
$\kappa$, let
\[ P_\calS=\setof{x_i}{i\in\kappa}\cup\setof{y_\alpha}{\alpha\in S}
\]\noindent
where $x_i$'s and $y_\alpha$'s are pairwise distinct, and let
$\leq_\calS$ be the \po\  on $P_\calS$ defined by
\[ \begin{array}{r@{}l}
p\leq_\calS q\  \Leftrightarrow\  {}&
                p=q\mbox{\  \  or }\\
                &p=x_i \mbox{ and }
                q=y_\alpha
                \mbox{ for some }i\in\kappa\mbox{ and }\alpha\in S
                \mbox{ with }i\in S_\alpha\  .
   \end{array}
\]\noindent
Let $B_\calS$ be the \Ba\  generated freely from $P_\calS$ except the
relation $\leq_\calS$. Note that the identity map on $P_\calS$
canonically induces a strictly order-preserving embedding of $P_\calS$
into $S_\calS$.

\begin{Prop}\label{nonwfn}
Suppose that $\cf(\kappa)\geq\omega_2$,
$S\subseteq\kappa$ is stationary \st\
$S\subseteq\setof{\alpha<\kappa}{\cf(\alpha)\geq\omega_1}$ and
$\calS=\setof{S_\alpha}{\alpha\in S}$ is \st\  $S_\alpha$ is a cofinal
subset of $\alpha$ for each $\alpha\in S$. If $P_\calS$ is embedded
into a \po\  $P$ by a strictly order-preserving mapping

then $P$ does not have the \wfn. In particular,
$B_\calS$ and its completion do not have the \wfn.
\end{Prop}
\prf
\Wolog, we may assume that
$P_\calS$ is a subordering of
$P$. Assume to the contrary that there is a \wfnm\
$\mapping{f}{P}{[P]^{\leq\aleph_0}}$. Let
\[
        \begin{array}{r@{}l}
                C=\setof{\xi<\kappa}{&\forall\eta<\xi\,\forall p\in 
F(x_\eta)\\
        &(\exists\alpha\in S\,(x_\eta\leq p\leq y_\alpha)\  \rightarrow\
        \exists\alpha'\in S\cap\xi\,(x_\eta\leq p\leq y_{\alpha'}))}.
        \end{array}
\]\noindent
Then $C$ is a club subset of $\kappa$. Let $\alpha\in C\cap S$ and let
\[ A=\setof{p\in F(y_\alpha)}{\exists\eta\in S_\alpha(p\in F(x_\eta)
                \  \land\  x_\eta\leq p\leq y_\alpha)}.
\]\noindent
Since $\alpha\in C$, for each $p\in A$ there is $\alpha_p<\alpha$ \st\
$p\leq y_{\alpha_p}$. Let
$\alpha^*=\sup\setof{\alpha_p}{p\in A}$. Since $A$ is countable we have
$\alpha^*<\alpha$. Let $\beta\in S_\alpha\setminus\alpha^*$. Since
$x_\beta\leq y_\alpha$, there is a $p\in A$ \st\
$x_\beta\leq p\leq y_\alpha$. Hence $x_\beta\leq y_{\alpha_p}$. But this
is impossible since $\alpha_p\leq \beta$.
\qedofProp

\section{Cohen models}\label{cohen-reals}
Consider the following principle:
\begin{assertion}{$({\ast}{}{\ast})$}
There is a sequence $\seqof{S_\alpha}{\alpha\in Lim(\omega_2)}$ \st\  each
$S_\alpha$ is a cofinal subset of $\alpha$ and that for any pairwise
disjoint $\seqof{x_\beta}{\beta<\omega_1}$ with
$x_\beta\in[\omega_2]^{<\aleph_0}$ for $\beta<\omega_1$, there are
$\beta_0<\beta_1<\omega_1$ \st\  $x_{\beta_0}\cap S_\alpha=\emptyset$ for
all $\alpha\in x_{\beta_1}\cap Lim(\omega_2)$ and that
$x_{\beta_1}\cap S_\alpha=\emptyset$ for
all $\alpha\in x_{\beta_0}\cap Lim(\omega_2)$ .
\end{assertion}

\begin{Prop}
Let $P=\Fn(\omega_2,2)$. Then 
$\forces{P}{({\ast}{}{\ast})}$.\label{starstar}
\end{Prop}
\prf \Wolog\  we may assume
$P=\Fn(\bigcup_{\alpha\in Lim(\omega_2)}\alpha\times\ssetof{\alpha}, 2)$.
For $\alpha\in Lim(\omega_2)$, let $\dotS_\alpha$ be a $P$-name \st\
$\forces{P}{\dotS_\alpha=\setof{\beta\in \alpha}{\dotg(\beta,\alpha)=1}}$
where $\dotg$ is the canonical name for the generic function. By 
genericity,
$\forces{P}{\dotS_\alpha\xmbox{ is cofinal in }\alpha}$ for every
$\alpha\in Lim(\omega_2)$.
Let $\dotS$ be
a $P$-name \st\
$\forces{P}{\dotS=\seqof{\dotS_\alpha}{\alpha\in Lim(\omega_1)}}$.

To show that $\dotS$ is forced to satisfy the property in 
\assertof{${\ast}{}{\ast}$},
let $\seqof{\dotx_\beta}{\beta<\omega_1}$ be a $P$-name of a sequence of
pairwise disjoint finite subsets of $\omega_2$. For each
$\beta<\omega_1$, let $p_\beta$ and $x_\beta\in[\omega_2]^{<\aleph_0}$ be
\st\  $p_\beta\forces{P}{\dotx_\beta=x_\beta}$. By thinning out the index
set $\omega_1$, we may assume \wolog\  that $\dom(p_\beta)$,
$\beta<\omega_1$ form a $\Delta$-system with the root $d$ and
$p_\beta\restr d$, $\beta<\omega_1$ are all equal to the same
$p\in P$. Since $p_\beta$, $\beta<\omega_1$ are then pairwise compatible,
$x_\beta$, $\beta<\omega_1$ are pairwise disjoint. Further, we may assume
also that $s_\beta$, $\beta<\omega_1$ form a $\Delta$-system with the
root $s$ where
$s_\beta=\setof{\gamma}{(\gamma,\alpha)\in\dom(p_\beta)\xmbox{ for some }
                \alpha<\omega_2}$.

Let $\beta_0<\beta_1<\omega_1$ be \st\
$x_{\beta_0}\cap s=\emptyset$, $x_{\beta_1}\cap s=\emptyset$,
$x_{\beta_0}\cap s_{\beta_1}=\emptyset$ and
$x_{\beta_1}\cap s_{\beta_0}=\emptyset$.
Let
\[
        \begin{array}{r@{}l}
        p^*={}&p_{\beta_0}\cup p_{\beta_1}
                \cup
                        \setof{((\beta,\alpha),0)}{\beta\in 
x_{\beta_0},\,\alpha
                        \in x_{\beta_1}\cap Lim(\omega_2)}\\[\jot]
                &\cup
                        \setof{((\beta,\alpha),0)}{\beta\in 
x_{\beta_1},\,\alpha
                        \in x_{\beta_0}\cap Lim(\omega_2)}
        \end{array}
\]\noindent
Then
$p^*\forces{P}{\dotx_{\beta_0}\cap\dotS_\alpha=\emptyset}$
for all $\alpha\in\dotx_{\beta_1}\cap Lim(\omega_2)$ and\
$p^*\forces{P}{\dotx_{\beta_1}\cap\dotS_\alpha=\emptyset}$
for all $\alpha\in\dotx_{\beta_0}\cap Lim(\omega_2)$.
\qedofProp
\begin{Prop}\mbox{}\label{ccc}
Suppose that $\seqof{S_\alpha}{\alpha\in Lim(\omega_2)}$ is as in
\assertof{${\ast}{}{\ast}$}. Let
$S=\setof{\alpha<\omega_2}{\cf(\alpha)=\omega_1}$ and
$\calS=\seqof{S_\alpha}{\alpha\in S}$. Then $B_\calS$ satisfies the c.c.c.
\end{Prop}
\prf
Otherwise we can find $I_\alpha\in[\omega_2]^{<\aleph_0}$,
$J_\alpha\in[S]^{<\aleph_0}$ for $\alpha<\omega_1$ and
$t(\alpha,i)$, $u(\alpha,\xi)\in\ssetof{+1,-1}$ for each
$i\in I_\alpha$, $\xi\in J_\alpha$ and $\alpha<\omega_1$ \st\
\[\textstyle
        z_\alpha=\prod_{i\in I_\alpha}t(\alpha,i)\,x_i
                \cdot\prod_{\xi\in J_\alpha}u(\alpha,\xi)\,y_\xi,
                \  \  \alpha<\omega_1\]\noindent
form a pairwise disjoint family of elements of
${B_\calS}^+$.

By $\Delta$-system argument, we may assume that
$I_\alpha\cup J_\alpha$, $\alpha<\omega_1$ are pairwise disjoint.
Applying \assertof{${\ast}{}{\ast}$} to
$\seqof{I_\alpha\cup J_\alpha}{\alpha<\omega_1}$, we find
$\beta_0<\beta_1<\omega_1$ \st\  $I_{\beta_0}\cap S_\xi=\emptyset$ for all
$\xi\in J_{\beta_1}$ and that $I_{\beta_1}\cap S_\xi=\emptyset$ for all
$\xi\in J_{\beta_0}$ . By definition of $B_\calS$, it follows that
$z_{\beta_0}\cdot z_{\beta_1}\not=0$\,.
This is a contradiction.
\qedofProp
\begin{Thm}\label{Cor1}
In a Cohen model (i.e.\  any model obtained by adding $\geq\aleph_2$
Cohen reals) there is a c.c.c.\  \cBa\  $B$ of density $\aleph_2$ without 
the \wfn.
\end{Thm}
\prf By \Propof{starstar}, \assertof{${\ast}{}{\ast}$} holds in a Cohen 
model.
Hence $B_\calS$ as in \Propof{ccc} satisfies the c.c.c. By \Propof{nonwfn},
the completion of $B_\calS$ does not have the \wfn.\qedofThm

\begin{Cor}\label{Cor2}
The \wfn\  of $(\psetof{\omega},{\subseteq})$ does not imply the
\wfn\  of all c.c.c.\  \cBas.
\end{Cor}
\prf If we start from a model of \CH\  and add $\aleph_2$ Cohen reals,
then $(\psetof{\omega},{\subseteq})$ has the \wfn\  in the resulting
model (see e.g.\  \cite{fuchino-soukup}). On the other hand, by 
\Thmof{Cor1}, there
is a c.c.c.\  \cBa\  without the \wfn\  in such a model. \qedofCor\qedskip

Under \CH, every c.c.c.\  \cBa\  of size $\aleph_2$ has the
\wfn\,(\cite{FuKoSh}). Hence it follows from the result above that \CH\
implies the negation of the principle \assertof{${\ast}{}{\ast}$}. This
can be also seen directly as follows:
\begin{Prop} \CH\  implies $\neg$\assertof{${\ast}{}{\ast}$}.
\end{Prop}
\prf Let $\seqof{S_\alpha}{\alpha\in Lim(\omega_2)}$ be any sequence \st\
each $S_\alpha$ is a cofinal subset of $\alpha$ for
$\alpha\in Lim(\omega_2)$. To show that
$\seqof{S_\alpha}{\alpha\in Lim(\omega_2)}$ is not as in
\assertof{${\ast}{}{\ast}$}, let $\chi$ be sufficiently large and let
$M\prec \calH(\chi)$ be \st\  $\cardof{M}=\aleph_1$;
$\seqof{S_\alpha}{\alpha\in Lim(\omega_2)}\in M$;
$\omega_1\subseteq M$; $\omega_2\cap M\in \omega_2$ and, letting
$\gamma=\omega_2\cap M$, $\cf(\gamma)=\omega_1$.
By \CH\  --- and since $\omega_1\subseteq M$ and $\cf(\gamma)=\omega_1$,
$[\gamma]^{\aleph_0}\subseteq M$.

Now choose by induction distinct $\alpha^0_\beta$,
$\alpha^1_\beta<\gamma$ for $\beta<\omega_1$ \st\
\assertof{1} $\alpha^0_\beta\in S_\gamma$, and
\assertof{2} $\setof{\alpha^0_\xi}{\xi<\beta}\subseteq S_{\alpha^1_\beta}$
for all $\beta<\omega_1$. \assertof{2} is possible: since
$\setof{\alpha^0_\xi}{\xi<\beta}\subseteq S_\gamma$ and
$\setof{\alpha^0_\xi}{\xi<\beta}\in M$, we have
\[M\models\exists\nu<\omega_2(\sup\setof{\alpha^1_\xi}{\xi<\beta}<\nu
        \land \setof{\alpha^0_\xi}{\xi<\beta}\subseteq S_\nu).\]
\noindent
Let $x_\beta=\ssetof{\alpha^0_\beta,\alpha^1_\beta}$ for
$\beta<\omega_1$. Then there are no $\beta_0<\beta_1<\omega_1$ \st\
$x_{\beta_0}\cap S_\alpha=\emptyset$ for all
$\alpha\in x_{\beta_1}\cap Lim(\omega_2)$. \qedofProp

\section{The Weak Freese-Nation property of c.c.c.\   complete Boolean
algebras under \GCH} \label{GCH}
In \cite{fuchino-soukup} it is proved that, assuming \CH\  and a weak
form of square principle at singular cardinals of cofinality $\omega$, 
every
c.c.c.\  \cBa\  has the \wfn. In this section we show that even \GCH\  does
not suffice for this result.

Hajnal, Juh\'asz and Shelah \cite{hajnal-juhasz-shelah} showed that,
starting from a model with a supercompact cardinal, a model of \GCH\  and 
the
following principle can be constructed:
\begin{assertion}{($\textstyle{{\ast}{}{\ast}{}{\ast}}$)}
There are a stationary
$S\subseteq\setof{\alpha<\omega_{\omega+1}}{\cf(\alpha)=\omega_1}$ and a
family $\calS=\seqof{S_\alpha}{\alpha\in S}$ \st\  each $S_\alpha$ is a 
cofinal
subset of $\alpha$ of ordertype $\omega_1$ and that, for all distinct
$\alpha$, $\beta\in S$, $S_\alpha\cap S_\beta$ is finite.
\end{assertion}
\begin{Prop}\mbox{}\label{ccc2}
Suppose that $\calS=\seqof{S_\alpha}{\alpha\in S}$ is as in
\assertof{${\ast}{}{\ast}{}{\ast}$}.
Then $B_\calS$ satisfies the c.c.c..
\end{Prop}
\prf
Otherwise we can find $I_\alpha\in[\omega_{\omega+1}]^{<\aleph_0}$,
$J_\alpha\in[S]^{<\aleph_0}$, $\alpha<\omega_1$ and $t(\alpha,i)$,
$u(\alpha,\xi)\in\ssetof{+1,-1}$ for each $\alpha<\omega_1$,
$i\in I_\alpha$ and $\xi\in J_\alpha$ \st
\[\textstyle z_\alpha=\prod_{i\in I_\alpha}t(\alpha,i)\,x_i
        \cdot\prod_{\xi\in J_\alpha}u(\alpha,\xi)\,y_\xi,\  \
\alpha<\omega_1
 \]\noindent
form a pairwise disjoint family of elements of
${B_\calS}^+$.

By $\Delta$-system argument, we may assume that
$I_\alpha\cup J_\alpha$, $\alpha<\omega_1$ are pairwise disjoint and each
$I_\alpha$ has the same size, say $n$.

For $\alpha<\beta<\omega_1$, since $z_\alpha\cdot z_\beta=0$, either
(I) there is $\eta\in J_\alpha$ \st\  $I_\beta\cap S_\eta\not=\emptyset$ or
else (II) there is $\xi\in J_\beta$ \st\
$I_\alpha\cap S_\xi\not=\emptyset$. If (I) holds then let us say that
$(\alpha,\beta)$ is of type (I).

Now, one of the following two cases should hold. We show that both of
them lead to a contradiction.

{\bf Case I. } There is an infinite subset $S$ of $\omega$ \st\  for every
$\beta\in\omega_1\setminus\omega$,
$\setof{k\in S}{(k,\beta)\xmbox{ is of type (I)}}$ is cofinite in $S$.
In this case, by thinning out the index set $\omega_1$, we may assume
that, for any $k\in\omega$ and $\beta\in\omega_1\setminus\omega$,
$(k,\beta)$ is of type (I). Since $\cardof{I_\alpha}=n$, for all
$\beta\in\omega_1\setminus\omega$, there are
$0\leq i^0(\beta)<i^1(\beta)<n+1$ \st\
$I^*_\beta=I_\beta\cap S_{i^0(\beta)}\cap S_{i^1(\beta)}\not=\emptyset$ by
Pigeonhole Principle.
Hence we can find an infinite $X\subseteq\omega_1\setminus\omega$ and
$0\leq i^0<i^1<n+1$ \st\  $i^0(\beta)=i^0$ and $i^1(\beta)=i^1$ for all
$\beta\in X$. But then
$\bigcup_{\beta\in X}I^*_\beta\subseteq S_{i^0}\cap S_{i^1}$. Since the
set on the left side is infinite as an infinite disjoint union of
non-empty sets, this is a contradiction to $({\ast}{\ast}{\ast})$.

{\bf Case II. } For any infinite subset $S\subseteq\omega$, there is
$\beta\in\omega_1\setminus\omega$ \st\  for infinitely many $k\in S$,
$(k,\beta)$ is not of type (I).
In this case, by thinning out the index set $\omega_1$, we may assume that
for each $k\in\omega$, there is
$\xi(k)\in J_\omega$ \st\  $I_k\cap S_{\xi(k)}\not=\emptyset$. Note that
$J_\omega$ is finite. So by thinning out further the first $\omega$ 
elements of the
index set $\omega_1$, we may assume that there is $\xi_0\in J_\omega$ 
\st\  $I_k\cap S_{\xi_0}\not=\emptyset$ for every $k<\omega$. Similarly we
may also assume that there are $\xi_i\in J_{\omega+i}$ for $1\leq i\leq n$ 
\st\
$I_k\cap S_{\xi_i}\not=\emptyset$ for every $k<\omega$.
For each $k<\omega$, we can find $i^0(k)<i^1(k)\leq n$ \st\
$I^*_k=I_k\cap S_{\xi_{i^0(k)}}\cap S_{\xi_{i^1(k)}}\not=\emptyset$ by
Pigeonhole Principle.
Since there are only $n(n-1)/2$ possibilities of
$i^0(k)<i^1(k)\leq n$, there are
$i^0<i^1\leq n$ and an infinite set $X\subseteq\omega$ \st\
for every $k\in X$, $i^0(k)=i^0$ and $i^1(k)=i^1$. It follows that
$S_{\xi_{i^0}}\cap S_{\xi_{i^1}}\supseteq\bigcup_{k\in X}I^*_k$. Since
$S_{\xi_{i^0}}\cap S_{\xi_{i^1}}$ is finite, this is a contradiction.
\qedofProp
\begin{Thm} It is consistent  with \GCH\  (modulo the consistency strength
of a \label{Cor3}supercompact cardinal) that there is a c.c.c.\
\cBa\  without the \wfn.
\end{Thm}
\prf Let $\calS$ be a family as in \assertof{${\ast}{}{\ast}{}{\ast}$}. By
\Propof{ccc2} and \Propof{nonwfn} the completion of $B_\calS$ is a
c.c.c.\  \cBa\  without the \wfn.\qedofThm\qedskip
\\
In \cite{fuchino-soukup} it is proved that under \CH\  and a very weak
version of the square principle at $\aleph_\omega$, every c.c.c.\  \cBa\ of
cardinality $\aleph_{\omega+1}$ has the \wfn. Hence we see that
consistency strength of some large cardinal is involved in
\assertof{${\ast}{}{\ast}{}{\ast}$}. This can be also seen directly as 
follows.

First let us review the weak form of the square principle used in
\cite{fuchino-soukup}.
$\Box^{*{}*{}*}_{{\aleph_1},\mu}$ is the following assertion: there exists 
a
sequence $\seqof{C_\alpha}{\alpha<\mu^+}$ and a club set
$D\subseteq\mu^+$ \st\  for $\alpha\in D$ with $\cf(\alpha)\geq{\omega_1}$
\begin{assertion}{}\mbox{}%
\lassert{y1} $C_\alpha\subseteq\alpha$, $C_\alpha$ is unbounded in
$\alpha$;\smallskip\\
\lassert{y2} 
$[\alpha]^{<{\omega_1}}\cap\setof{C_{\alpha'}}{\alpha'<\alpha}$
dominates $[C_\alpha]^{<{\omega_1}}$ (\wrt\  $\subseteq$).
\end{assertion}
It can be easily seen that
$\Box^{*{}*{}*}_{{\aleph_1},\mu}$ follows from the very weak square 
principle for
$\mu$ by Foreman and Magidor \cite{foreman-magidor} (see
\cite{fuchino-soukup}).
\begin{Prop}\label{astastast}
$2^{\aleph_0}<\aleph_\omega$ and
$\Box^{*{}*{}*}_{{\aleph_1},\omega_\omega}$ implies the negation of
\assertof{${\ast}{}{\ast}{}{\ast}$}.
\end{Prop}
\prf Let $\seqof{C_\alpha}{\alpha<{\omega_\omega}+}$ and
$D\subseteq{\omega_\omega}^+$ be as in the definition of
$\Box^{*{}*{}*}_{{\aleph_1},\omega_\omega}$. Suppose that $S$ and
$\seqof{S_\alpha}{\alpha\in S}$ are as in 
\assertof{${\ast}{}{\ast}{}{\ast}$}.
We show that there are $\xi$, $\eta\in S$, $\xi\not=\eta$ \st\
$S_\xi\cap S_\eta$ is not finite. By replacing $S$ by $S\cap D$ and
$C_\alpha$ by $C_\alpha\cap S_\alpha$ for $\alpha\in S\cap D$, we may
assume that $S\subseteq D$ and $\seqof{S_\alpha}{\alpha\in S}$ is just
the sequence $\seqof{C_\alpha}{\alpha\in S}$.

By \assertof{y2}, for each $\alpha\in S$, there is
$\beta_\alpha<\alpha$ of countable cofinality \st\
$C_\alpha\cap C_{\beta_\alpha}$ is infinite. By Fodor's theorem we may
assume that every $\beta_\alpha$, $\alpha\in S$ are the same $\beta$.
Since there are only $<\aleph_\omega$ different subsets of the countable
set $C_\beta$, we can find $\xi$, $\eta\in S$, $\xi\not=\eta$ \st\
$S_\xi\cap C_\beta= S_\eta\cap C_\beta$. It follows that
$S_\xi\cap S_\eta$ is infinite. \qedofProp

\section{The \wfn\  of c.c.c.\  \cBas\  under weak square
principles}\label{square}
In this section, we investigate the c.c.c.\  \cBas\  $B$ for which we can 
prove
(in \ZFC\  or in some extension of it) that the \wfn\
of $B$ is equivalent with the \wfn\  of $(\psof{\omega},{\subseteq})$.
\Lemmaof{equivalence} was an easy observation in this direction. S.\
Koppelberg observed in \cite{koppelberg} that the Cohen algebra
$\bbd{C}(\omega)$ is such a \Ba.

Since $(\psof{\omega},{\subseteq})$ can be embedded in every \cBa, it
follows from \Propof{criterions} \assertof{a} and one of the remarks
before it that, if $(\psof{\omega},{\subseteq})$ does not have the \wfn\
then no \cBa\  can have the \wfn.

\begin{Prop} The measure algebra $\bbd{R}(\omega)$ has the \wfn\  if and
\label{measure-alg}
only if $(\psof{\omega},{\subseteq})$ has the \wfn.
\end{Prop}
\prf By the remark above and by \Propof{criterions},\assertof{f} together
with one of the remarks before \Propof{criterions}, it is
enough to find a strictly order-preserving embedding of $\bbd{R}(\omega)$ 
into
$(\fnsp{\omega}{\Re},{\leq})$.  We may replace $\omega$ by the countable 
set
$Clop(\fnsp{\omega}{2})$ where $Clop(\fnsp{\omega}{2})$
denotes the clopen sets of the Cantor space $\fnsp{\omega}{2}$.

We define
$\mapping{e}{\bbd{R}(\omega)}{\fnsp{Clop(\fnsp{\omega}{2})}{\Re}}$ by
taking
$e(a)(c)={\mu}(a\cap c)
$
for $c\in Clop(\fnsp{\omega}{2})$. Then clearly $e$ is a
order-preserving map of $\bbd{R}(\omega)$ into
$\fnsp{Clop(\fnsp{\omega}{2})}{\Re}$.

To show that $e$ is a strictly order-preserving, assume that
${\mu}(a\setminus b)>0$.
Then there is $c\in Clop(\fnsp{\omega}{2})$ such that
${\mu}((a\setminus b)\cap c)>{\mu}(c)/2$. Then
$e(a)(c)> {\mu}(c)/2> e(b)(c)$, so $e(a)\not\le e(b)$.
\qedofProp\qedskip
\\
In a similar way, we can also show that $\Re(\omega)$ is a retract of
$\psof{\omega}$ as a partially ordered set though it is known that
$\Re(\omega)$ is {\em not} a retract of $\psof{\omega}$ as a \Ba: the
mapping
$\mapping{e}{\Re(\omega)}{\psof{Clop(\fnsp{\omega}{2})\times\rationals}}$ 
defined by
$e(c)=\setof{(a,q)}{a\in Clop(\fnsp{\omega}{2}),\,q<\mu(a\cap c)}$ is
easily seen to be a strictly order preserving embedding.

In general, the \wfn\  of $\bbd{C}(\kappa)$ or that of $\bbd{R}(\kappa)$
for arbitrary $\kappa\geq\omega$ is not equivalent with the \wfn\
of $(\psof{\omega},\subseteq)$. In the next section we shall
give a model where $(\psof{\omega},\subseteq)$ has the \wfn\  while
$\bbd{C}(\aleph_\omega)$ (and hence also $\bbd{R}(\aleph_\omega)$) does 
not.

However the equivalence does hold if $\kappa<\aleph_\omega$ or some
consequences of $\neg0^\#$ are available.
To prove this, we need
the following instance of Theorem 7 in \cite{fuchino-soukup}:
\begin{Thm}{\rm (Theorem 7 in \cite{fuchino-soukup}
for $\kappa=\aleph_1$) }\label{fuchino-soukup}
Suppose that $\mu>\cf(\mu)=\omega$,
$\cf([\lambda]^{\aleph_0},\subseteq)=\lambda$ for cofinally many 
$\lambda<\mu$ and
$\Box^{{\ast}{\ast}{\ast}}_{\aleph_1,\mu}$ holds. Then for any
sufficiently large regular $\chi$ and $x\in\calH(\chi)$, there is a
matrix $(M_{\alpha,i})_{\alpha<\mu^+,i<\omega}$ \st
\begin{itemize}
\item[{\rm(1)}] $M_{\alpha,i}\prec\calH(\chi)$, $x\in M_{\alpha,\beta}$,
$\omega_1\subseteq M_{\alpha,i}$ and $\cardof{M_{\alpha,i}}<\mu$ for all
$\alpha<\mu^+$ and $i<\omega$;
\item[{\rm(2)}] $(M_{\alpha,i})_{i<\omega}$ is an increasing sequence for 
each
$\alpha<\mu^+$;
\item[{\rm(3)}] If $\alpha<\mu^+$ and $\cf(\alpha)\geq\omega_1$, then there 
is an
$i^*<\omega$ \st\  for every $i^*\leq i<\omega$,
$[M_{\alpha,i}]^{\aleph_0}\cap M_{\alpha,i}$ is cofinal in
$([M_{\alpha,i}]^{\aleph_0},{\subseteq})$;
\item[{\rm(4)}] Let $M_\alpha=\bigcup_{i<\omega}M_{\alpha,i}$ for
$\alpha<\mu^+$. Then $M_\alpha\prec\calH(\chi)$ (by \assertof{1} and
\assertof{2}). Moreover $(M_\alpha)_{\alpha<\mu^+}$ is continuously
increasing and $\mu^+\subseteq \bigcup_{\alpha<\mu^+}M_\alpha$. \qed
\end{itemize}\end{Thm}

For a \cBa\  $B$ and $X\subseteq B$, let us denote by
$\langle X\rangle^{\rm cm}_B$ the
complete subalgebra of $B$ generated completely by $X$.
\begin{Thm}
Let $\lambda$ be a cardinal. \label{ccccBa}
Suppose that for every $\mu<\lambda$ with $\mu>\cf(\mu)=\omega$, we have:
\smallskip
\\
\assert{i} $\cf([\mu]^{\aleph_0},\subseteq)=\mu^+$;\\
\assert{ii} $\Box^{{\ast}{\ast}{\ast}}_{\aleph_1,\mu}$\,.\smallskip
\\
Then for any c.c.c.\  \cBa\  $B$ with a complete generator of size
$<\lambda$, $B$ has the \wfn\  if and only if every complete subalgebra of
$B$ generated completely by a countable subset of $B$ has the \wfn.
\end{Thm}
\prf ``Only if'' part of the theorem follows from \Propof{criterions}
\assertof{b}. ``If'' part of the theorem is proved by induction on the
minimal cardinality of a subset $X$ of $B$ completely generating $B$.

If $X$ is countable, then there is nothing to prove. Let
$\mu=\cardof{X}<\lambda$ and suppose that we have the theorem for any
c.c.c.\  \cBa\  with a complete generator of size $<\mu$.

If $\cf(\mu)>\omega$, then letting $X=\setof{x_\alpha}{\alpha<\mu}$ and
$B_\beta=\langle\setof{x_\alpha}{\alpha<\beta}\rangle^{\rm cm}_B$ for
$\beta<\mu$, we
have $B_\beta\leq_{c}B$ and hence $B_\beta\leq_\sigma B$ for all
$\beta<\mu$. By induction hypothesis, every $B_\beta$, $\beta<\mu$ has
the \wfn. By the c.c.c.\  of $B$ $B_\gamma=\bigcup_{\beta<\gamma}B_\beta$ 
for
all limit $\gamma<\mu$ with $\cf(\gamma)>\omega$ and also
$B=\bigcup_{\beta<\mu}B_\beta$. Hence by \Propof{criterions} \assertof{c},
it follows that $B$ has the \wfn.

If $\cf(\mu)=\omega$, then there is
$(M_{\alpha,i})_{\alpha<\mu^+,i<\omega}$ as in the previous theorem for
$x=(B,X)$.

For $\alpha<\mu^+$ and $i<\omega$, let
$B_{\alpha,i}=\langle B\cap M_{\alpha,i}\rangle^{\rm cm}_B$ and
$B_\alpha=\bigcup_{i<\omega}B_{\alpha,i}$.
\begin{Claim}
For every $\alpha<\mu^+$, $B_\alpha$ has the \wfn\  and
$B_\alpha\leq_\sigma B$.
\end{Claim}
\prfofClaim For every $i<\omega$, $B_{\alpha,i}$ has the \wfn\  by
induction hypothesis. Since $B_{\alpha,i}\leq_c B$ for every
$i<\omega$ it follows that $B_\alpha\leq_\sigma B$. Also by
\Propof{criterions} \assertof{c} it follows that $B_\alpha$ has the \wfn.
\qedofClaim
\begin{Claim}
If $\gamma<\mu^+$ and $\cf(\gamma)>\omega$, then
$B_\gamma=\bigcup_{\alpha<\gamma}B_\gamma$.
\end{Claim}
\prfofClaim Suppose that $a\in B_\gamma$. Then, by the c.c.c.\  of $B$,
there is an $i<\omega$ and $s\in[B\cap M_{\gamma,i}]^{\aleph_0}$ \st\
$a\in\langle s\rangle^{\rm cm}_B$. By \assertof{3} in
\Thmof{fuchino-soukup}, we may
assume that $s\in M_{\gamma,i}$. By \assertof{4}, there is
$\alpha<\gamma$ and $j<\omega$ \st\  $s\in M_{\alpha, j}$. It follows that
$s\subseteq M_{\alpha,j}$ and hence $a\in B_{\alpha,j}\subseteq B_\alpha$.
\qedofClaim
\begin{Claim}
$B=\bigcup_{\alpha<\mu^+}B_\alpha$\,.
\end{Claim}
\prfofClaim
By the last statement of
\assertof{4} in \Thmof{fuchino-soukup} and \assertof{i},
$[X]^{\aleph_0}\cap \bigcup_{\alpha<\mu^+}M_\alpha$ is cofinal in
$([X]^{\aleph_0},{\subseteq})$.

Suppose now that $a\in B$. Then by the c.c.c.\  of $B$, there is a 
countable
$s\in[X]^{\aleph_0}$ \st\  $a\in\langle s\rangle^{\rm cm}_B$. By the
remark above, we
may assume that $s\in \bigcup_{\alpha<\mu^+}M_\alpha$, say
$s\in M_{\alpha^*,i^*}$ for some $\alpha^*<\mu^+$ and $i^*<\omega$. Then
$s\subseteq B\cap M_{\alpha^*,i^*}$ and hence
$a\in B_{\alpha^*,i^*}\subseteq B_{\alpha^*}$.\qedofClaim\smallskip

Now by \Thmof{criterions} \assertof{c} and the claims above, it follows
that $B$ has the \wfn.
\qedofThm
\begin{Cor}\label{equivalence2}
Let $\lambda$ be as in \Thmof{ccccBa} and $\omega\leq\kappa<\lambda$. Then
\tfae:\smallskip
\\
\assert{0} $(\psof{\omega},{\subseteq})$ has the \wfn;
\\
\assert{1} $\bbd{C}(\kappa)$ has the \wfn;
\\
\assert{2} $\bbd{R}(\kappa)$ has the \wfn.
\end{Cor}
\prf \assertof{1} \Implies \assertof{0} and \assertof{2}
\Implies \assertof{0} follows from \Propof{criterions} \assertof{b}.

For \assertof{0} \Implies \assertof{2} assume that
$(\psof{\omega},{\subseteq})$ has the \wfn. Since $\bbd{R}(\kappa)$ has
the complete generator $Clop(\fnsp{\kappa}{2})$, it is
enough by \Thmof{ccccBa} to show that every subalgebra of
$\bbd{R}(\kappa)$ has the \wfn. Let $A$ be such an algebra then there is
$X\subseteq[\kappa]^{\aleph_0}$ \st\  $A\leq_c \bbd{R}(X)$. Since
$\bbd{R}(X)$ has the \wfn\  by \Propof{measure-alg}, it follows by
\Propof{criterions} \assertof{b} that $A$ also has the \wfn.\qedofCor
\qedskip

Note that the conditions on $\lambda$ in \Thmof{ccccBa} hold vacuously for
$\lambda=\aleph_\omega$. Hence we obtain the following as a special case of 
the
corollary above:
\begin{Cor}
\Tfae\  (in \ZFC):\label{ccccBa-ZFC}
\smallskip
\\
\assert{0} $(\psof{\omega},{\subseteq})$ has the \wfn;
\\
\assert{1} $\bbd{C}(\aleph_n)$ has the \wfn\  for some/all $n\in\omega$;
\\
\assert{2} $\bbd{R}(\aleph_n)$ has the \wfn\  for some/all 
$n\in\omega$.\qed
\end{Cor}
\section{Chang's Conjecture}\mbox{}\label{chang}
In this section, we give a negative answer to Problem 3 mentioned in the
introduction (see \Thmof{chang-pomega} below) and show that
\Corof{ccccBa-ZFC} in the last section is the optimal assertion among what 
we can
obtain in \ZFC.

\begin{Thm}\mbox{}\label{chang-pomega}%
Suppose that $V_0$ is a transitive model of \ZFC\  \st\
\[ V_0\models\mbox{\rm\GCH\  +\  }
        (\aleph_{\omega+1},\aleph_\omega)\pfeil(\aleph_1,\aleph_0)\,.
\]\noindent
Let $P$ be a c.c.c.\  \po\  in $V_0$ of cardinality $\aleph_1$ adding a
dominating real. Let
$\eta\in\fnsp{\omega}{\omega}$ be a dominating real over $V_0$
generically added by $P$ and let $V_1=V_0[\eta]$. Note that 
$V_1\models\GCH$.
In $V_1$ let $Q=\Fn(\aleph_\omega,\omega)$ and let $\dot Q$ be a
corresponding $P$-name. Then we have:
\[ V_1\models\forces{Q}{(\psof{\omega},{\subseteq})
                \mbox{ does not have the \wfn}}.
\]\noindent
\end{Thm}
\prf In the proofs of this and next theorems, we shall denote by a dotted
symbol a name of an element in a generic extension. By the same symbol
without the dot, we denote the corresponding element in a fixed generic
extension. Without further mention, we shall identify $P*\dotQ$ names
with the corresponding $Q$-name in $V_1$ and vice versa.

Now toward a contradiction, assume that there is a $Q$-name $\dotF$ in 
$V_1$ \st\
\begin{markedformula}{($\otimes$)}
V_1\models\forces{Q}{%
        \dotF\mbox{ is a \wfnm\  over }(\fnsp{\omega}{\omega},{\leq^*})}.
\end{markedformula}
Let $\dotvarphi$ be a $P*\dotQ$-name of the function
$\aleph_\omega\,\rightarrow\,\omega$
generically added by $Q$ over $V_1$.
Let $V_2=V_1[\varphi]$. By \GCH, we can find in $V_0$ a scale
$\seqof{f_\alpha}{\alpha<\aleph_{\omega+1}}$ in
$(\prod_{n\in\omega}\aleph_n,{\leq^*})$. \Wolog, we may assume that for
every $\alpha<\aleph_{\omega+1}$ and $n\in\omega$,
$f_\alpha(n)\in\aleph_n\setminus\aleph_{n-1}$ (where we set
$\aleph_{-1}=0$). For each $\alpha\in\aleph_{\omega+1}$, let
\[ g_\alpha=\varphi\circ f_\alpha\,.
\]\noindent

Let $\chi$ be sufficiently large and let $N\prec(\calH(\chi),{\in})$ be 
\st\
$N$ contains every thing we need in the course of the proof,
$\cardof{\aleph_\omega\cap N}=\aleph_0$ and
$otp(\aleph_{\omega+1}\cap N)=\omega_1$ --- the latter two conditions are
possible by $(\aleph_{\omega+1},\aleph_\omega)\pfeil(\aleph_1,\aleph_0)$.

In $V_0$, let $\setof{\xi_{n,k}}{k\in\omega}$ be an enumeration of
$(\aleph_n\setminus\aleph_{n-1})\cap N$ for each $n\in\omega$. Here again, 
we
use the convention that $\aleph_{-1}=0$. Let $\doth^*$ be a
$P*\dotQ$-name of an element of $\fnsp{\omega}{\omega}$ \st
\[ 
\forces{P*\dotQ}{\doth^*(n)=\max\setof{\dotvarphi(\xi_{n,k})}{k<\doteta(n)}
        \mbox{ for all }n\in\omega}.
\]\noindent
\begin{Claim}\label{claim1}
For every $\alpha\in\aleph_{\omega+1}\cap N$,
$\forces{P*\dotQ}{\dotg_\alpha\leq^*\doth^*}$.
\end{Claim}
\prfofClaim
Since $\alpha\in N$ we have $f_\alpha\in N$. Hence there is a function
$\mapping{e_\alpha}{\omega}{\omega}$ in $V_0$ \st\
$f_\alpha(n)=\xi_{n,e_\alpha(n)}$\,. Since $\eta$ is dominating, there is
$n^*\in\omega$ \st
\[ V_1\models e_\alpha\restr\omega\setminus n^*\leq 
\eta\restr\omega\setminus n^*.
\]\noindent
By definition of $\doth^*$, it follows that
\[ V_2\models g_\alpha(n)=\varphi\circ f_\alpha(n)
                =\varphi(\xi_{n,e_\alpha(n)})\leq h^*(n)
\]\noindent
for all $n\geq n^*$.
\qedofClaim\qedskip

Let $N_0=N$, $N_1=N_0[\eta]$ and $N_2=N_1[\varphi]$. Then we
have $N_2\prec \calH(\chi)[\eta][\varphi]$.

Let $\doth_n\in N_0$, $n\in\omega$ be $P*\dotQ$-names \st
\[ \forces{P*\dotQ}{\setof{\doth_n}{n\in\omega}=\dotF(\doth^*)\cap\dotN_2}.
\]\noindent
Let $S_n\in[\aleph_\omega]^{\aleph_0}\cap N_0$ be \st, regarding
$\doth_n$ as a $P$-name of $\dotQ$-name,
\[ V_1\models\forces{P}{\doth_n\mbox{ is a }\Fn(S_n,\omega)\mbox{-name}}.
\]\noindent
This is possible since $P$ has the c.c.c.\  and
$\forces{P}{\dotQ\mbox{ has the c.c.c.}}$.
For each $n\in\omega$, let $s_n\in\prod_{n\in\omega}\aleph_n\cap N_0$ be
defined (in $N_0$) by $s_n(k)= \sup S_n\cap\aleph_k$ for $k\in\omega$.
Since $\seqof{f_\alpha}{\alpha<\aleph_{\omega+1}}$ was taken to be a
scale on $\prod_{n\in\omega}\aleph_n$, for each $n\in\omega$ there is
$\alpha_n\in\aleph_{\omega+1}\cap N_0$ \st\  $s_n\leq^*f_{\alpha_n}$. Let
$\alpha^*\in\aleph_{\omega+1}\cap N_0$ be \st\
$\sup\setof{\alpha_n}{n\in\omega}\leq\alpha^*$.

Now, by the choice of $\doth_n$, $n\in\omega$, the following claim
contradicts to ($\otimes$) and \Claimof{claim1}, and hence proves the
theorem:
\begin{Claim}
$V_1\models\forces{Q}{\dotg_\alpha\,{\not\leq}^*\doth_n\mbox{ for all 
}n\in\omega}$.
\end{Claim}
\prfofClaim Assume to the contrary that, in $V_1$, we have
\[ q\forces{Q}{\dotg_\alpha\restr(\omega\setminus k)\leq
                \doth_n\restr(\omega\setminus k)}
\]\noindent
for some $q\in Q$ and $n$, $k\in\omega$. We may assume that
\[ s_n\restr\omega\setminus k<f_{\alpha^*}\restr\omega\setminus k
\]\noindent
and $\sup(\dom(p))\leq f_\alpha(m)$ for all $m\in\omega\setminus k$ as
well. Working further in $V_1$,
let $m^*=k+1$ and let $q'\leq q$ be \st\
\[ q'\forces{Q}{\doth_n(m^*)=j^*}
\]\noindent
for some $j^*\in\omega$. Then $q''=q\cup(q'\restr S_n)$ also forces the 
same
statement. Since
$f_{\alpha^*}(m^*)\in\aleph_{m^*}\setminus(s_n(m^*)\cup\dom(p))$, we have
\[ f^*_\alpha(m^*)\not\in \dom(q'').
\]\noindent
Hence
\[ q^*=q''\cup\ssetof{(f_{\alpha^*}(m^*),j^*+1)}
\]\noindent
is an element of $Q$ and $q^*\leq q''\leq q$. But
\[ q^*\forces{Q}{\dotg_{\alpha^*}(m^*)
        =\dotvarphi\circ f_{\alpha^*}(m^*)
        =j^*+1>j^*=\doth_n(m^*).
}
\]\noindent
This is a contradiction.\qedofClaim\\
\qedofThm
\begin{Thm}\mbox{}\label{cohen-chang}
Suppose that $V_0$, $P$, $\eta$, $\doteta$, $V_1$ are as in the previous
theorem. Then
\[ V_1\models{\bbd C}(\aleph_\omega)\mbox{ does not have the \wfn.}
\]\noindent
\end{Thm}
\prf
Suppose to the contrary that $F\in V_1$ is a \wfnm\  on
${\bbd C}(\aleph_\omega)$. Let $X=\setof{x_\xi}{\xi<\aleph_\omega}$ be a
free subset of ${\bbd C}(\aleph_\omega)$ completely generating the whole
${\bbd C}(\aleph_\omega)$. We may take $X\in V_0$.

Let
$\seqof{f_\alpha}{\alpha<\aleph_{\omega+1}}$ be as in the proof of the
previous theorem. For $n\in\omega$ let $c_n=x_n\cdot-\sum_{m<n}x_m$.
$\setof{c_n}{n\in\omega}$ is then a partition of ${\bbd C}(\aleph_\omega)$.
For
$\alpha<\aleph_{\omega+1}$ let
\[ b_{\alpha,n}={\sum}_{m>n}(x_{f_\alpha(m)}\cdot c_m)
\]\noindent
Let $\chi$ be sufficiently large and let $N\prec(\calH(\chi),{\in})$ be 
\st\
$N$ contains every thing we need in the course of the proof,
$\cardof{\aleph_\omega\cap N}=\aleph_0$ and
$otp(\aleph_{\omega+1}\cap N)=\omega_1$. In $V_0$, let
$\setof{\xi_{n,k}}{k\in\omega}$ be an enumeration of
$(\aleph_n\setminus\aleph_{n-1})\cap N$ for each $n\in\omega$.
In $N_1=N[\eta]$ let
\[ b^*={\sum}_{n\in\omega}(c_n\cdot {\sum}_{l<\eta(n)}x_{\xi_{n,l}})\,.
\]\noindent
\begin{Claim}
For every $\alpha\in\aleph_{\omega+1}\cap N$, there is
$n_\alpha\in\omega$ \st\  $b_{\alpha,n_\alpha}\leq b^*$.
\end{Claim}
\prfofClaim
In $V_0$, let $s_\alpha\in\fnsp{\omega}{\omega}$ be \st\
$f_\alpha(n)=\xi_{n,s_\alpha(n)}$ for all $n\in\omega$. Since $\eta$ is
dominating, there is an $n_\alpha\in\omega$ \st\
$s_\alpha\restr\omega\setminus n_\alpha\leq\eta\restr\omega\setminus 
n_\alpha$. By
definition of $b_{\alpha,n}$ it follows that $b_{\alpha,n_\alpha}\leq b^*$.
\qedofClaim\qedskip

Now, let $\seqof{d_n}{n\in\omega}$ be an enumeration of
$F(b^*)\cap\bbd{C}(\aleph_\omega)\restr b^*\cap N$ and $\dotd_n$,
$n\in\omega$ be corresponding $P$-names. We
can choose these names so that $\seqof{\dotd_n}{n\in\omega}\subseteq N$.

For each $n$, there is $S_n\in[\aleph_{\omega+1}]^{\aleph_0}\cap N$ \st\
\[ \forces{P}{\dotd_n\in\langle S_n\rangle^{\rm 
cm}_{\bbd{C}(\aleph_\omega)}}.
\]\noindent
Let $f^*_n\in\prod_{n\in\omega}\aleph_n\cap N$ be defined by
\[ f^*_n(m)=\min(\aleph_m\setminus S_n)
\]\noindent
for $m\in\omega$. Since $\seqof{f_\alpha}{\alpha<\aleph_{\omega+1}}$ is
a scale on $(\prod_{n\in\omega}\aleph_n,{\leq^*})$, there is
$\alpha^*_n\in\aleph_{\omega+1}\cap N$ \st\
$f^*_n\leq^*f_{\alpha^*_n}$. Let $\alpha^*\in\aleph_{\omega+1}\cap N$ be
\st\  $\alpha^*_n<\alpha^*$ for all $n\in\omega$.

Then similarly to the proof of the previous theorem. The following claim 
gives
a desired contradiction:
\begin{Claim}
For every $n\in\omega$, $b_{\alpha^*,n_{\alpha^*}}\not\leq d_n$.
\end{Claim}
\prfofClaim For $n\in\omega$, let $m\in\omega$ be \st\
$n_{\alpha^*}<m$, $c_m\cdot -d_n\not=0$ and
$f^*_n(m)<f_{\alpha^*}(m)$. Then $f_{\alpha^*}(m)\not\in S_n$.
Hence
$x_{f_{\alpha^*}(m)}\cdot c_m\not\leq d_n$ but
$x_{f_{\alpha^*}(m)}\cdot c_m\leq b_{\alpha^*,n_{\alpha^*}}$ by
definition of $b_{\alpha^*,n_{\alpha^*}}$. \qedofClaim\\
\qedofThm

\mbox{}\vfill\mbox{}\hfill
\parbox[t]{10cm}{
\noindent
{\bf Authors' addresses}\bigskip\bigskip\\
{\rm Saka\'e Fuchino}
\smallskip\\
Kitami Institute of Technology, Kitami, Japan.\\
{\tt fuchino@math.cs.kitami-it.ac.jp}\medskip\\
{\rm Stefan Geschke}\smallskip\\
II. Mathematisches Institut, Freie Universit\"at Berlin, Germany. \\
{\tt geschke@math.cs.kitami-it.ac.jp}\medskip\\
{\rm Saharon Shelah}\smallskip\\
Institute of Mathematics, The Hebrew University of Jerusalem, 91904
Jerusalem, Israel.\\
Department of Mathematics, Rutgers University, New Brunswick, NJ 08854,
U.S.A.\\
{\tt shelah@math.huji.ac.il}\medskip\\
{\rm Lajos Soukup}\smallskip\\
Kitami Institute of Technology, Kitami, Japan.\\
Institute of Mathematics, Hungarian Academy of Science, Budapest, 
Hungary.\\
{\tt soukup@math.cs.kitami-it.ac.jp}}


\begin{thebibliography}{99}
\newcommand{\bysame}{\underline{\  \  \  \  \  \  }}
\newcommand{\vol}{Vol.\  }
\newcommand{\no}{No.\  }
\iftesting
\let\Bibitem\bibitem
\def\bibitem#1{\Bibitem{#1}\marginpar{{\tiny #1}}}
\fi
\bibitem{blass}
A.\  Blass, Combinatorial cardinal characteristics of the continuum, to
appear in {\it Handbook of Set-Theory}.
\bibitem{brendle-fuchino-soukup}
J.\  Brendle, S.\  Fuchino, L.\  Soukup: Coloring ordinals by reals, 
preprint.
\bibitem{foreman-magidor}
M.\  Foreman and M.\  Magidor, A very weak square principle,.
\bibitem{FuKoSh}
S.\  Fuchino, S.\  Koppelberg and S.\  Shelah, Partial orderings with the
weak Freese-Nation property, Annals of Pure and Applied Logic, 80 (1996),
35--54.
\bibitem{fuchino-soukup}
S.\  Fuchino and L.\  Soukup,
More set-theory around the weak Freese-Nation property,
Fundamenta Matematicae 154 (1997), 159--176.
\bibitem{fuchino-geschke-soukup}
S.\  Fuchino, S.\  Geschke and L.\  Soukup, On the weak Freese-Nation
property of $\psof{\omega}$, submitted.
\bibitem{geschke} S.\  Geschke, On $\sigma$-filtered \Bas, Dissertation,
Berlin (1999).
\bibitem{heindorf-shapiro}
L.\  Heindorf and L.B.\  Shapiro,
{\it Nearly Projective Boolean Algebras}, Springer Lecture Notes of
Mathematics 1596, (1994).
\bibitem{hajnal-juhasz-shelah}
A.\  Hajnal, I.\  Juh\'asz and S.\  Shelah, Splitting strongly almost 
disjoint
families,
Transactions of the American Mathematical Society 295 No.1 (1986), 369-397.

\bibitem{koppelberg}
S.\  Koppelberg, Applications of $\sigma$-filtered Boolean algebras,
Advances in Algebra and Model Theory, eds. M. Droste, R. G"obel,
Gordon and Breach, sc. publ. (1997), 199-213.
\bibitem{kunen}
K.~Kunen: {\it Set Theory}, North-Holland (1980).


\end{thebibliography}
\end{document}